\documentclass[runningheads]{llncs}
\usepackage{graphicx} 

\usepackage{amsfonts}
\usepackage{amsmath}

\usepackage{amsthm}

\usepackage{authblk}
\usepackage{xcolor}
\usepackage{trimclip}
\usepackage{float}
\newcommand{\hsplit}[3][\par]{%
\clipbox{0em 0em #2 0em}{#3}#1
\clipbox{{\the\dimexpr\width-#2\relax} 0em 0em 0em}{#3}}

\usepackage[T1]{fontenc}
\usepackage{graphicx}
\begin{document}
\title{Generalization of Ramanujan's Famous Nested Radicals to the nth Root and their Evaluation}
\titlerunning{Nested Radicals}
\author{Narendranath Mitra\orcidID{0009-0007-9056-7658}}
\authorrunning{Narendranath Mitra}
\institute{Jadavpur University, Kolkata 700032, West Bengal, India}
\maketitle              
\begin{abstract}
\fontsize{9}{12}\selectfont \sloppy SRINIVASA RAMANUJAN posed a problem on infinite nested radical of the square root in the Journal of Indian Mathematical Society in 1911. He had generated the problem years before in the form of an example illustrating a more general theorem. Of course, how we'd figure the solution out without Ramanujan's theorem was scarcely obvious. Generation of infinite nested radicals, the roots being of higher order i.e. proceeding from the cube roots to the nth root as well as their evaluation may seem to be somewhat complex. For this purpose a general solution in the form of a theorem needs to be provided whereby symmetric patterns of infinite nested radicals of any order of roots may be evaluated including that of Ramanujan involving square root. In our work, we have rendered the required solution. 

\end{abstract}
\section{General Solution to Nested Radicals}
The title indicates an endeavour to define a general solution to the nested radical involving nth root following Ramanujan's famous identity. Legendary Indian mathematician Srinivasa Ramanujan's first mathematics paper was published in the Journal of Indian Mathematical Society in 1911. In his article (\cite{Ramanujan}, p. 323, Question 289) he proposed two problems involving infinite nested radicals of the square root. One of them was to find the value of:

\[ \sqrt{1+2\sqrt{1+3\sqrt{1+4\sqrt{1+...}}}} \]

Ramanujan himself offered its solution - it was simply 3. This problem was generated from the following general form which Ramanujan (\cite{Berndt}, p. 108, Entry 4) had discovered.
Let $x, n'$ and $a$ denote arbitrary numbers. Then
\begin{equation} \label{eq1}
x+n'+a = \sqrt{ax + (n'+a)^2 + x\sqrt{a(x+n') + (n'+a)^2 + (x+n')\sqrt{...}}} 
\end{equation}
Putting $ x=2, n'=1 $ and $a=0$ in (\ref{eq1}) yields that
\begin{equation} \label{eq2}
\sqrt{1+2\sqrt{1+3\sqrt{1+4\sqrt{1+...}}}} = 3
\end{equation}

Titu Andreescu, Răzvan Gelca (\cite{Andree}, p. 112) and Christopher G. Small (\cite{Small}, p. 21-24 \& 72-75) approached the problem considering the function $ f: {[1,\infty)} \mapsto \mathbb{R} $,
\begin{equation} \label{eq3}
f(x) = \sqrt{1+x\sqrt{1+(x+1)\sqrt{1+(x+2)\sqrt{1+...}}}}
\end{equation}
being a nested radical involving square roots and proved that $ f(x)= x+1 $; the particular case $x=2$ provides Ramanujan's identity (\ref{eq2}).

Next, we observe the following:

Let us consider the function $ f: {[1,\infty)} \mapsto \mathbb{R} $ be defined as 
\begin{equation} \label{eq4}
f(x)=x+1
\end{equation}
\subsection{Nested square root}
Squaring both sides of (\ref{eq4}),
\[ (f(x))^2 = 1+2x+x^2 = 1+x(x+2) = 1+xf(x+1) \]
i.e. \centerline{ $ f(x) = \sqrt{ 1+xf(x+1) } 
= \sqrt{1+x\sqrt{1+(x+1)f(x+2)}} $ }
\[ = \sqrt{1+x\sqrt{1+(x+1)\sqrt{1+(x+2)f(x+3)}}} = ... \] 
i.e. \centerline{ $ f(x) = \sqrt{1+x\sqrt{1+(x+1)\sqrt{1+(x+2)\sqrt{...}}}} $ }
Then
\begin{equation} \label{eq5}
x+1 = \sqrt{1+x\sqrt{1+(x+1)\sqrt{1+(x+2)\sqrt{...}}}}
\end{equation}
\subsection{Nested cube root}
Cubing both sides of (\ref{eq4}),
\[ (f(x))^3 = 1+3x+3x^2+x^3 = 1+3x+x^2(x+3) = 1+3x+x^2f(x+2) \]
i.e. \centerline{ $ f(x) = \sqrt[3]{1+3x+x^2f(x+2)} $}
\[ = \sqrt[3]{1+3x+x^2\sqrt[3]{1+3(x+2)+(x+2)^2f(x+4)}} = ... \] 
i.e. \centerline{ $ f(x) = \sqrt[3]{1+3x+x^2\sqrt[3]{1+3(x+2)+(x+2)^2\sqrt[3]{...}}} $ }
Then
\begin{equation} \label{eq6}
x+1 = \sqrt[3]{1+3x+x^2\sqrt[3]{1+3(x+2)+(x+2)^2\sqrt[3]{...}}}
\end{equation}
\subsection{Nested fourth root}
Raising both sides of (\ref{eq4}) to fourth power, 
\[ (f(x))^4 = 1+4x+6x^2+4x^3+x^4 = 1+4x+6x^2+x^3(x+4) = 1+4x+6x^2+x^3f(x+3) \]
i.e. \centerline{ $ f(x) = \sqrt[4]{1+4x+6x^2+x^3f(x+3)} $}
\[ = \sqrt[4]{1+4x+6x^2+x^3\sqrt[4]{1+4(x+3)+6(x+3)^2+(x+3)^3f(x+6)}} = ... \] 
i.e. \centerline{ $ f(x)= \sqrt[4]{1+4x+6x^2+x^3\sqrt[4]{1+4(x+3)+6(x+3)^2+(x+3)^3\sqrt[4]{...}}} $ }
Then
\begin{equation} \label{eq7}
x+1= \sqrt[4]{1+4x+6x^2+x^3\sqrt[4]{1+4(x+3)+6(x+3)^2+(x+3)^3\sqrt[4]{...}}}
\end{equation}
\subsection{Nested nth root}
Raising to nth power, and proceeding similarly,
\begin{equation*} 
\begin{split}
(f(x))^n & = 1+\binom{n}{1}x+\binom{n}{2}x^2+...+\binom{n}{n-2}x^{n-2}+\binom{n}{n-1}x^{n-1}+\binom{n}{n}x^n \\
& = 1+\binom{n}{1}x+\binom{n}{2}x^2+...+\binom{n}{n-2}x^{n-2}+x^{n-1}f(x+n-1)
\end{split}
\end{equation*}
\begin{flushleft}
Now we generalize this problem for solving the nth root nested radical:
\end{flushleft}

\centerline{ $ f(x) = \sqrt[n]{1+\binom{n}{1}x+\binom{n}{2}x^2+...+\binom{n}{n-2}x^{n-2}+x^{n-1}f(x+n-1)} $}
\begin{flushleft}
Then
\end{flushleft}

\hsplit[\par\hspace{9.5em}]{22.3em}{ $ \displaystyle\ f(x) = \sqrt[n]{1+\binom{n}{1}x+...+\binom{n}{n-2}x^{n-2}+x^{n-1}\sqrt[n]{1+\binom{n}{1}(x+n-1)+...+\binom{n}{n-2}(x+n-1)^{n-2} +(x+n-1)^{n-1}\sqrt[n]{...}}} $ }
\begin{flushleft}
Now we achieve the evaluation of the nested radicals of the nth root with the following theorem.
\end{flushleft}

\fbox{
\begin{minipage}{34em}
\newtheorem*{theorem*}{Theorem}
\theoremstyle{plain}
\begin{theorem*}
Let x be a non-negative real number in ${[1,\infty)}$. Then

\hsplit[\par\hspace{9.2em}]{21.9em} { $ \displaystyle\ x+1 = \sqrt[n]{1+\binom{n}{1}x+...+\binom{n}{n-2}x^{n-2}+x^{n-1}\sqrt[n]{1+\binom{n}{1}(x+n-1)+...+\binom{n}{n-2}(x+n-1)^{n-2} +(x+n-1)^{n-1}\sqrt[n]{...}}} $ }

\end{theorem*}
\end{minipage}
}

\begin{flushleft}
 Putting $n=2$ and $x=2$, the above theorem readily yields Ramanujan's problem (\ref{eq2}).  
\end{flushleft}
\section{Polynomial Nature of Nested Cube Root}

Let us justify the equation (\ref{eq6}) involving the nested radical containing cube roots.

Let us consider the function $ f: {[1,\infty)} \mapsto \mathbb{R} $ be defined as

\begin{equation} \label{eq9}
f(x) = \sqrt[3]{1+3x+x^2\sqrt[3]{1+3(x+2)+(x+2)^2\sqrt[3]{...}}}
\end{equation}
\begin{equation} \label{eq10}
\text{i.e.}
f(x)= \sqrt[3]{1+3x+x^2\sqrt[3]{7+3x+(x+2)^2\sqrt[3]{13+3x+(x+4)^2\sqrt[3]{...}}}} 
\end{equation}
Now it is required to prove that $ f(x)=x+1 $.

\subsection{Using Algebra}
We verify this using some sort of Algebra.

Cubing (\ref{eq9}) we obtain the functional equation:
\begin{equation} \label{eq11}
(f(x))^3 = 1+3x+x^2f(x+2) 
\end{equation}

If $ f(x) $ be a polynomial of degree $ n $, then L.H.S. of our functional equation (\ref{eq11}) would be of degree $ 3n $ and the R.H.S. of degree $ n+2 $, whence
\[ 3n = n+2  \quad \textrm{giving} \quad n = 1 , \]
providing us that the polynomial $ f(x ) $ is of degree one.

So, let $ f(x) = ax+b ; $ therefore, by (\ref{eq11}),
\[ a^3x^3 + 3a^2bx^2 + 3ab^2x + b^3 = ax^3 + (2a+b)x^2 + 3x +1 \]

Equating coefficients of like powers of $ x $,
\[ a^3=a, 3a^2b = 2a+b, 3ab^2 = 3, b^3 = 1 , \]
and it is evident that $ a = 1 $ and $ b = 1 $ satisfy these four relations.

\subsection{Using Analysis}

Now we consider the following function:
\[ f(x)= \sqrt[3]{1+3x+x^2\sqrt[3]{7+3x+(x+2)^2\sqrt[3]{13+3x+(x+4)^2\sqrt[3]{...}}}} \tag{9} \] 

Truncating to cube roots in progression, we obtain a monotonically increasing sequence of nested radicals. All we require to show is that this sequence is bounded from above and below whereby the matter of convergence is settled.
\[f(x) \leq \sqrt[3]{(1+3x+x^2)\sqrt[3]{(7+3x+(x+2)^2)\sqrt[3]{(13+3x+(x+4)^2)\sqrt[3]{...}}}} \]
\[ = \sqrt[3]{(1+3x+x^2)\sqrt[3]{(11+7x+x^2)\sqrt[3]{(29+11x+x^2)\sqrt[3]{...}}}} \]
\[ \leq \sqrt[3]{(x+1)(x+2)\sqrt[3]{(x+3)(x+4)\sqrt[3]{(x+5)(x+6)\sqrt[3]{...}}}}\]
\[ \leq \sqrt[3]{1(x+1).2(x+1)\sqrt[3]{3(x+1).4(x+1)\sqrt[3]{5(x+1).6(x+1)\sqrt[3]{...}}}}\]
\[ = (x+1)\sqrt[3]{1.2\sqrt[3]{3.4\sqrt[3]{5.6\sqrt[3]{...}}}} 
\leq (x+1)\sqrt[3]{4\sqrt[3]{16\sqrt[3]{36\sqrt[3]{...}}}} \]
\[ \leq (x+1)\sqrt[3]{4\sqrt[3]{16\sqrt[3]{64\sqrt[3]{...}}}} 
= (x+1).4^{3/4} 
= 2^{3/2}(x+1)\]
also, \[ f(x) \geq \sqrt[3]{x^2\sqrt[3]{x^2\sqrt{x^2\sqrt[3]{...}}}} = x \geq 2^{-3/2}(x+1) \] 

Combining both, we obtain
\begin{equation} \label{eq12}
2^{-3/2}(x+1) \leq f(x) \leq 2^{3/2}(x+1)
\end{equation}

And replacing $ x $ by $ x+2 $, we get 
\[ 2^{-3/2}(x+3) \leq f(x+2) \leq 2^{3/2}(x+3) \]

Hence from the preceding double inequality, and using the functional relation (\ref{eq11})
\[ (f(x))^3 = 1+3x+x^2f(x+2) \]
it yields
\[ 1+3x+2^{-3/2}x^2(x+3) \leq (f(x))^3 \leq 1+3x+2^{3/2}x^2(x+3) \]
which implies
\[ 2^{-3/2}(1+3x)+2^{-3/2}x^2(x+3) \leq (f(x))^3 \leq 2^{3/2}(1+3x)+2^{3/2}x^2(x+3) \]
i.e. \centerline{ $ 2^{-3/2}(x+1)^3 \leq (f(x))^3 \leq 2^{3/2}(x+1)^3 $ } 
Then 
\begin{equation} \label{eq13}
2^{-1/2}(x+1) \leq f(x) \leq 2^{1/2}(x+1)
\end{equation}

Compare (\ref{eq13}) with (\ref{eq12}). Repeating successively the argument, we find that
\[ (2^{-1/2})^{1/3^k}(x+1) \leq f(x) \leq (2^{1/2})^{1/3^k}(x+1) \]

Letting $ k\to\infty $ gives us
\[ x+1 \leq f(x) \leq x+1 \]
And hence \centerline{$ f(x) = x+1 $,}
which proves the result.

\subsection{General Solution and Special Case}
Our general solution for the cube root nested radical is:
\[ x+1 = \sqrt[3]{1+3x+x^2\sqrt[3]{7+3x+(x+2)^2\sqrt[3]{13+3x+(x+4)^2\sqrt[3]{...}}}} \tag{9} \]

Symmetry can be seen everywhere in this equation. We have three basic terms under each radical. The first term has the following progression: 1, 7, 13 etc. 

The second term, $ 3x $, does not change while third term progresses as: $ x^2, (x+2)^2, (x+4)^2 $ etc.

Putting $ x= 1 $ in our identity (\ref{eq10}), we obtain a special case,
\begin{equation} \label{eq14}
\sqrt[3]{4+1^2\sqrt[3]{10+3^2\sqrt[3]{16+5^2\sqrt[3]{22+...}}}} = 2
\end{equation}

\begin{figure}[H]
\renewcommand{\figurename}{Figure}
\renewcommand{\thefigure}{}
    \centering
    \includegraphics[width=10cm]{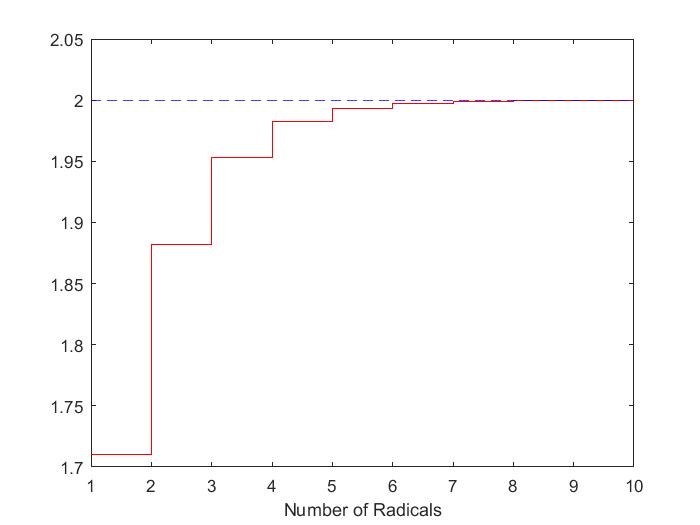}
    \caption{Graphical representation of the special case }
    \label{fig:radical}
\end{figure}
\begin{flushleft}
The figure shows that the above sequence of nested radicals converges to 2.    
\end{flushleft}

\noindent \textit{Notes and Comments.} Finding a value to nested square root that initially seemed so mysterious was solved by Ramanujan. It would be still amazing to evaluate a nested cube root and in generality a nested nth root. Here we have established a general solution to the nth root nested radicals following Ramanujan's pattern.

\bibliographystyle{plain}
\bibliography{name}

\end{document}